\documentclass{amsart}
\usepackage{amsmath,amsthm,amssymb,xypic,array}

\theoremstyle{plain}  
\newtheorem{Th}{Theorem}
\newtheorem{theorem}{Theorem}[section]

\newtheorem{proposition}[theorem]{Proposition}
\newtheorem{lemma}[theorem]{Lemma}   
\newtheorem{corollary}[theorem]{Corollary}

\theoremstyle{definition}

\newtheorem{definition}[theorem]{Definition}

\theoremstyle{remark}

\newtheorem{claim}{Claim}

\newtheorem{remark}[theorem]{Remark}

\DeclareMathOperator{\Bsl}{Bs}

\DeclareMathOperator{\mult}{mult}

\DeclareMathOperator{\fr}{Frac}

\newcommand{\QED}{\ifhmode\unskip\nobreak\fi\quad {\rm Q.E.D.}} 

\newcommand\frup[1]{\lceil{#1}\rceil}
\newcommand\frdown[1]{\lfloor{#1}\rfloor}

\newcommand\kdn[3]{\lceil\frac{{#1\choose #2}}{#3}\rceil}
\newcommand\bin[2]{{#1\choose #2}}
\newcommand\Span[1]{\langle{#1}\rangle}

\newcommand\iso{\cong}

\newcommand{\G}{\mathcal{G}}

\renewcommand{\H}{\mathcal{H}} 

\newcommand{\I}{\mathcal{I}}

\renewcommand{\O}{\mathcal{O}} 
\renewcommand{\P}{\mathbb{P}}

\newcommand{\rat}{\dasharrow}

\newcommand{\exf}{{expected and effective}}
\newcommand{\win}{\sl {\bf wIn}}
\newcommand{\nbhd}{{neighborhood}\ }

\begin{document}

\title{Base loci of linear systems and the Waring Problem}

\author{
Massimiliano Mella}
\address{
Dipartimento di Matematica\\ Universit\`a di Ferrara\\44100 Ferrara
 Italia}
\email{mll@unife.it}

\date{September 2007}
\thanks{Partially supported by Progetto PRIN 2006 ``Geometria
  sulle variet\`a algebriche'' MUR}
\maketitle

\section*{Introduction}
The Waring problem for forms is the quest for an additive decomposition
of homogeneous polynomials into powers of linear ones. The subject has
been widely considered in old times, \cite{Sy}, \cite{Hi}, \cite{Ri}
and \cite{Pa}, with special regards to the existence of a unique decomposition of
this type. The following, see \cite{RS}, was the
state of the art at the beginning of ${\rm XX^{th}}$ century.
A general form $f$ of degree $d$ in $n + 1$ variables has a unique
presentation as a sum of 
$s$  powers of linear forms in the following cases:
\begin{itemize}
\item $n = 1$, $d = 2k - 1$ and $s = k$, \cite{Sy}
\item $ n = 3$, $d = 3$ and $s = 5$ Sylvester's Pentahedral Theorem
  \cite{Sy}
\item $n = 2$, $d = 5$ and $s = 7$ \cite{Hi}, \cite{Ri}, \cite{Pa}.
\end{itemize}

After these remarkable results the Waring problem had kind of a
rest and only quite recently come back on the main scene, see
\cite{Ci} for a complete account of all contributions ranging from
algebraic geometry to commutative algebra and computer science. 

I do expect that the one listed are the only possible cases in which
the decomposition is unique.
A partial result in this direction is the following.
\begin{Th}[\cite{Me}]Let $f$ be a general homogeneous form of degree $d$ in
  $n+1$ variables. Assume that $d>n>1$. Then $f$ is
  expressible as sum of $d$th powers of linear forms in a
  unique way if and only if 
  $n=2$ and $d=5$.
\label{th:1}
\end{Th}

The aim of this note is to give more evidence to the above expectation.

\begin{Th}\label{th:main} Let $f$ be a general homogeneous form of degree $d$ in
  $n+1$ variables. Assume that $d\geq 5$, $n\geq 3$, and
  $\frac{\bin{n+d}{n-1}}{n}$ is an integer. Then $f$ is
  never expressible as sum of $d$th powers of linear forms in a
  unique way.
\end{Th}

Let me briefly comment upon the two theorems. The idea is to translate
the original statement into one on birational maps of $\P^n$. To have
theorem \ref{th:1} it is then enough to study the singularities of
special linear systems $\G_{d,n,l}$, with prescribed double points. This is done by a
degeneration argument essentially borrowed from \cite{AH1}. To go
beyond, namely $d\leq n$, one has to study the base locus of these
special linear systems. The main difficulty is that the degeneration
technique seems to be hopeless in this realm. To study the base locus
of a linear system of projective dimension $n$ on an $n$-fold, one has
to take track of all the elements of a base. While, in general, the
degeneration  technique
allows you to better understand a sublinear system discarding
completely the other elements. Instead of trying to determine the base
locus of $\G_{d,n,l}$ I prove a weaker statement on the base locus on
an open Zariski set containing the imposed points. This, under the
divisibility assumption, is then enough to conclude.

The paper is organised as follows. I first introduce the main notations
and preliminaries. Then I study the base
locus of these special linear systems and bound the degree of the maps
to prove Theorem \ref{th:main}.

\section{Notations and preliminaries}
\label{sec:not}

Unless otherwise stated I work over the field of complex
numbers. First I introduce what is needed to study
linear systems with prescribed singularities.
\begin{definition} Let $p\in\P^n$ be a point. The double point at $p$
  in $\P^n$ is the scheme given by the square of the ideal sheaf of
  $p$. If $P\subset\P^n$ is a collection of points, I denote by $P^2$
  the double points supported on $P$. In particular the linear system 
$|\I_{P^2}(d)|$ is given by hypersurfaces of degree $d$
  singular at $P$.
\end{definition}

Given a collection of points $X=Q\cup Q_H$, with $Q_H$ supported
on a hyperplane $H$, let $\tilde{X}$ be the residual of $X^2$ with respect
to $H$. 
That is $\tilde{X}=Q^2\cup Q_H$. Then there is the Castelnuovo exact
sequence given by
$$0\to\I_{\tilde{X}}(d-1)\to\I_{X^2}(d)\to\I_{Q^2,H}(d)\to 0$$
This gives the following sequence on cohomology
\begin{equation}
\label{eq:exseq}
0\to H^0(\P^n,\I_{\tilde{X}}(d-1))\to H^0(\P^n,\I_{X^2}(d))\to H^0(H,\I_{Q_H^2}(d))
\end{equation}

\begin{definition} Consider a collection $P$, of $l$ general points in
  $\P^n$.
Define 
$$\G_{d,n,l}:=|\I_{P^2}(d)|$$
Fix a hyperplane $H$ and a collection $X=Q\cup Q_H$, where $Q$ is given by $(l-h)$ general points in
$\P^n$, and $Q_H=\cup_1^h q_j$ is given by $h$ general points in $H$.
Define
$$\H_{H,d,n,l,h}:=|\I_{X^2}(d)|$$
\label{def:h}

\end{definition}

In this paper I am interested in non empty linear systems of type
$\G_{d,n,l}$, for this I introduce the following definition.

\begin{definition}
 I say that the linear system $\G_{d,n,l}$ is  {\sl expected} if
$$\dim\G_{d,n,l}=\bin{n+d}{n}-(n+1)l-1$$
Moreover if $\G_{d,n,l}$ is expected and $\dim\G_{d,n,l}\geq 0$ I say
that it is 
{\sl expected and effective} 

Note that  if
$\G_{d,n,l}$ is {\exf} then $\G_{d,n,l'}$ is {\exf} for any $l'<l$.
Similarly for linear systems of type $\H$. I say that $\H_{H,d,n,l,h}$ is  {\sl expected and
    effective} 
 if 
$$\dim\H_{H,d,n,l,h}=\bin{n+d}{n}-(n+1)l-1\geq 0$$
Note that if $\H_{H,d,n,l,h}$ is {\exf} then by semi-continuity
$\G_{d,n,l}$ is {\exf}.
In the following I frequently ask {\exf} linear systems of type $\H$ to
satisfy the following further properties. The linear system 
$\H_{H,d,n,l,h}$ is what I need ({\win}) if 
\begin{itemize}
\item $\H_{H,d,n,l,h}$ is {\exf},
\item $|\H_{H,d,n,l,h}\otimes \I_H|\neq \emptyset$
\item   (\ref{eq:exseq}) is a short exact sequence.
\end{itemize}
\end{definition}

I am interested in studying the base loci of linear systems of
type $\G_{d,n,l}$. To do this I use a degeneration argument.

\begin{lemma}[\cite{Me}]\label{le:dimsing} Let $\Delta$ be a complex disk around
  the origin. 
Consider the product $V=\P^n\times \Delta$,
 with the natural projections, $\pi_1$ and $\pi_2$. Let $V_t=\P^n\times\{t\}$ and
 $\O_V(d)=\pi_1^*(\O_{\P^n}(d))$.
Fix a configuration $q_1,\ldots,q_l$ of $l$ points on $V_0$ and let
 $\sigma_i:\Delta\to V$ be sections such that $\sigma_i(0)=q_i$ and
 $\{\sigma_i(t)\}_{i=1,\ldots,l}$ are general points of $V_t$ for
$t\neq 0$. Let $P_t=\cup_{i=1}^l\sigma_i(t)$.

 Consider the linear system $\H_t=|\O_{V_t}(d)\otimes\I_{P_t^2}|$. 
Assume that $l<\kdn{d+n}{n}{n+1}$ and $\dim\H_0=\dim\H_t$, for $t\in\Delta$. Let
$\psi_i(t):=\dim_{\sigma_i(t)}\Bsl\H_t$. 
Then for $t\neq 0$,  we have 
$$\psi_i(t)\leq \min\{j|\psi_j(0)\}$$
\end{lemma}

The way to pass from linear system of type $\H$ to those of type $\G$
is the content of the following.

\begin{corollary}\label{co:degbsl} Assume that
  $\H_{H,d,n,l,h}$ is {\win} and $\Bsl\H_{H,d,n,l,h}=q_i^2$ in a
  neighbourhood of a point $q_i$. Then $\Bsl(\G_{d,n,l})=P^2$ in a
  neighbourhood of $P$.
\end{corollary}

\begin{proof} I am in the condition to apply Lemma
  \ref{le:dimsing}, thus $\psi_i(t)=0$ for any $i$. Moreover the blow
  up of $q_i$ solves $\Bsl\H_{H,d,n,l,h}$ in a neighbourhood of
  $q_i$. This means that, with
  Lemma's notations,the blow up of
  $\sigma_i(\Delta)$ resolves $\Bsl\G_{d,n,l}$ in a neighbourhood of $p_i$. A monodromy
  argument shows that this is true for any point $p_i\in P$.
\end{proof}

I recall a way to translate the existence of a unique decomposition into a
statement on the map defined by linear system of type $\G_{d,n,k}$.

\begin{proposition}[\cite{Me}] Let $f$ be a general homogeneous form of degree $d$ in
  $n+1$ variables, with $n\geq 2$.  Then $f$ is
  expressible as a sum of $(k+1)$ powers of linear forms in a
  unique way  only if the map associated to the linear system
  $\G_{d,n,k}$ is birational, where $k=\frac{\bin{d+n}{n}}{n+1}-1$ is an integer.
\label{pro:bridge}
\end{proposition}

To test this condition, when $d\leq n$, I have to understand the base locus of the
linear system $\G_{d,n,k}$. In the following I shall do it by a
degeneration argument. For this let me recall the following results
from \cite{Me}.

\begin{lemma} Fix a hyperplane $H\subset\P^n$. Let 
$$l_d:=l=\kdn{n+d+1}{n}{n+1}-\kdn{n+d}{n-1}{n}$$
and 
$$h_d:=h=\kdn{n+d+1}{n}{n+1}-\kdn{n+d}{n-1}{n}-\kdn{n+d}{n}{n+1}+\kdn{n-1+d}{n-1}{n}$$ 
Assume that $d\geq 4$ and $n\geq 3$ 
then $\H_{H,d,n,l,h}$ is {\win}.
\label{le:lh}
\end{lemma}

\section{Base loci of linear systems}

First I study the base locus of $\G_{d,n,l}$ when $l$ is quite small.

\begin{theorem} \label{th:bsl} Assume that $d\geq 4$, $n\geq 3$, and
  $l\leq \kdn{n+d+1}{n}{n+1}-\kdn{n+d}{n-1}{n}$.
 Then $\G_{d,n,l}$ is {\exf} and $\Bsl(\G_{d,n,l})=P^2$ in a
 neighbourhood of $P$.
\end{theorem}
\begin{proof}
It is clear that it is enough to prove it for $l=\kdn{n+d+1}{n}{n+1}-\kdn{n+d}{n-1}{n}$.
I prove the claim by  induction on  $n$. 

Let $l'=\kdn{n+d}{n}{n+1}-\kdn{n+d-1}{n-1}{n}$ and 
$h=l-l'$.
 To be able to apply induction
on $n$ I need to prove that $\gamma=h-
\kdn{n+d}{n-1}{n}-\kdn{n-1+d}{n-2}{n-1}\leq 0$.
Note that 
$$
h=\kdn{n+d+1}{n}{n+1}-\kdn{n+d}{n-1}{n}-\kdn{n+d}{n}{n+1}+\kdn{n+d-1}{n-1}{n}<
\frac{\bin{n-1+d}{n-1}}n-\frac{\bin{n+d}{n-1}}{n(n+1)}+2$$
so that
\begin{eqnarray*}\gamma=h-(\kdn{n+d}{n-1}{n}-\kdn{n-1+d}{n-2}{n-1})<\frac{\bin{n+d-1}{n-1}}n+\frac{\bin{n+d-1}{n-2}}{n-1}-\bin{n+d}{n-1}\frac{n+2}{n(n+1)}+3\\
<\frac{(n+d-1)\cdot\ldots\cdot(d+2)}{n!}(1-\frac{n+d}{n+1})+3=\alpha
\end{eqnarray*}
Standard computations give $\alpha\leq 0$ in the following cases
\begin{itemize}
\item[.] $d\geq 4$ and $n\geq 9$
\item[.] $d\geq 5$ and $n\geq 6$
\item[.] $d\geq 6$ and $n\geq 4$
\item[.] $d\geq 9$ and $n\geq 3$
\end{itemize}
In the remaining cases I have to compute $\gamma$ directly
$${\setlength{\extrarowheight}{4pt}
\begin{tabular}{|c|cccccc|}
\hline
$(d,n)$&(4,3)&(4,4)&(4,5)&(4,6)&(4,7)&(4,8)\\
\hline
$\gamma$&-1&0&-3&-1&-3&-2\\
\hline
$(d,n)$&(5,3)&(5,4)&(5,5)&(6,3)&(7,3)&(8,3)\\
\hline
$\gamma$&-2&-2&-2&-1&-1&-5\\
\hline
\end{tabular}}
$$

This shows that I can apply induction hypothesis on 
$\G_{d,n-1,h}$.
 Fix a general hyperplane $H\subset\P^n$, then 
by Lemma \ref{le:lh} the linear system $\H_{H,d,n,l,h}$ is {\win}. In
particular this proves that
$\G_{d,n,l}$ and $\G_{d-1,n,l-h}$ are {\exf}.

 Let $Q$ be $l-h$ general points in $\P^n$. The ,linear system $\G_{d-1,n,l-h}$ is
 expected and effective and $\dim\G_{d-1,n,l-h}\geq n+1$. 
Let $D_1,\ldots,
D_{n-1}\in\G_{d-1,n,l-h}$ be general divisors and $H$ an
hyperplane. Let $C\subset D_1\cap\ldots\cap D_{n-1}$ be a general one dimensional
irreducible component. The linear system $\H_{H,d,n,l,h}$ is ${\win}$,
 therefore I can choose the hyperplane $H$ in such a way
that
\begin{itemize}
\item[.] $X=Q\cup Q_H$ and $Q_H\subset C\cap H$
\item[.] $H$ is tangent to $C$ at a point, say
$q_1\in Q_H$.
\end{itemize}
 Let $\epsilon:Y\to \P^n$ be the blow up of $q_1$ with
exceptional divisor $E$,
$H_Y=\epsilon^{-1}_*H$, and $\H_Y=\epsilon^{-1}_*\H_{H,d,n,l,h}$. The
choice of $H$ yields
$$\Bsl\H_Y\cap E\subset H_Y\cap E$$
 The linear
system $\H_{H,d,n,l,h}$  is {\win}, in particular all elements of
 $\G_{d,n-1,h}$ lift to elements in $\H_{H,d,n,l,h}$. By
induction hypotheses the linear system $\G_{d,n-1,h}$
has the correct base locus in a neighbourhood of $q_1$. Then
 $\Bsl(\H_{H,d,n,l,h})=q_1^2$ in a neighbourhood of $q_1$.
I then conclude by Corollary \ref{co:degbsl}.

It is left the proof of the first step of induction. That is the
statement for $n=3$. Arguing as in  the induction step it is enough to prove that 
$\G_{d,2,h}$ has the expected base locus. Moreover 
$$h<\frac{d+2}3+\frac32<\bin{\frdown{\frac{d}2}+2}{2}-1$$
 Then there are reducible divisors
$D=D_1\cup D_2\in \G_{d,2,l}$ such that  $P\in D_i$. This means that
the base locus in a neighbourhood of $P$ is the one prescribed. 
\end{proof}

\begin{remark} Note that an important feature of the above proof is that the
linear system $\H_{H,d,n,l,h}$ breaks in two parts $\G_{d-1,n,l-h}$
and $\G_{d,n-1,h}$ both of dimension at least $n+1$. This is of course
not true  when one imposes the maximal number of double points. The
divisibility condition in Theorem \ref{th:main} force all the linear
system ``on one side''.
\end{remark}

\begin{theorem} Assume that $d+1\geq 5$, $l=\frac{\bin{n+d+1}{n}}{n+1}-1$ and
  $\frac{\bin{n+d}{n-1}}{n}$ are integers. Then
  $\Bsl(\G_{d+1,n,l})=P^2$ in a neighbourhood of $P$. Assume that
  $n\geq 3$ then the map given by the
  linear system $\G_{d+1,n,l}$ is not birational.
\label{th:fcbsl}
\end{theorem}
\begin{proof}
Let
 $h=\frac{\bin{n+d}{n-1}}{n}$. In \cite[Proof of Theorem 4.3]{Me} it
 is pr oven that  $\H_{H,d+1,n,l,h}$ is
{\win}. Note that the numerical hypothesis yields $\dim
 \G_{d,n,l-h}=\dim\H_{H,d+1,n,l,h}= n+1$. That is
 $\H_{H,d+1,n,l,h}\subseteq H+\G_{d,n,l-h}$. Let $P$ be a length $l-h$
 0-scheme and $X=P\cup Q_H$ a length $l$ 0-scheme.
By Theorem \ref{th:bsl}
$\Bsl(\G_{d,n,l-h})=P^2$ in a \nbhd of $P$.
 Let $D_1,\ldots D_{n-1}\in \G_{d,n,l-h}$
be general elements and $C= D_1\cap\ldots\cap D_{n-1}$ the irreducible
 curve passing through $P$. 
Let $H$ be an hyperplane, I can assume that $Q_H\subset C\cap H$. The
 general choice of $H$ allows me to assume that 
 the general element in $|\G_{d,n,l-h}\otimes\I_{Q_H}|$ is not tangent to
 $C$ at $p_1$. This shows that $\Bsl(\H_{H,d+1,n,l,h})=p_1^2$ in a
 neighbourhood of $p_1$. Then Corollary \ref{co:degbsl} yields
 $\Bsl(\G_{d+1,n,l})=P^2$  in a neighbourhood of $P$.

 Let $\chi:\P^n\rat\P^n$ the map given by the linear system
 $\H_{H,d+1,n,l,h}$ or equivalently
 $|\G_{d,n,l-h}\otimes \I_{Q_H}|$.
I will decompose $\chi$ in several steps.
Let $\epsilon:Z\to \P^n$ be the resolution of $\Bsl\G_{d,n,l-h}$, with
$\G_Z=\epsilon^{-1}_*\G_{d,n,l-h}$ and $H_Z=\epsilon^{-1}_*H$. 
I have the following:
\begin{itemize}
\item[.] $E_1,\ldots,E_{l-h}$ the exceptional divisor over $P$ are all
  projective spaces of dimension $n-1$,
\item[.] the curve $C_Z=\epsilon^{-1}_* C$ is a curve section of $G_Z$.
\end{itemize}
Let $\Psi:Z\to\P^{n+h}$ be the morphism associated to the linear
system $\G_Z$.

  \[
 \xymatrix{
   &Z\ar[dl]_{\epsilon}\ar[r]^{\Psi}&\P^{n+h} \\
 \P^n&              &}
 \]

Then $\Psi(E_i)$ is a, eventually projected, 2-Veronese embedding of
$\P^{n-1}$. This bounds the  dimension of the span of $\Psi(E_i)$
\begin{equation}
\label{eq:span}\dim\Span{\Psi(E_i)}<\bin{n+1}2\end{equation}
I already proved that $\Bsl(\G_{d,n,l-h})=p_1^2$  in a neighbourhood of
$p_1$. Then I have
\begin{equation}
\label{eq:birat}
\#(\Psi(C_Z\cap E_1))=\G_Z^{n-1}\cdot E_1=2^{n-1}
\end{equation}
and by construction $\Span{\Psi(C_Z)}=\P^{h+1}$.

To complete the map  $\chi$ I have to project from a
linear space $\Lambda\iso\P^{h-1}$ $h$-secant to $\Psi(H_Z)$.

Assume that the map $\chi$ is birational. Then the projection
from $\Lambda$, say $\Pi_{\Lambda}$, is
birational. Let $\Lambda=\Span{\Lambda_1,x}$ for some point $x$ and
linear space $\Lambda_1$. Let me further factor this projection
with the projection from $\Lambda_1$ and then from $x$

  \[
 \xymatrix{
   &Z\ar[dl]_{\epsilon}\ar[r]^{\Psi}&\P^{n+h}\ar@{.>}[d]_{\Pi_{\Lambda}}\ar@{.>}[r]^{\Pi_{\Lambda_1}}&
   \P^{n+1}\ar@{.>}[dl]^{\Pi_{x}}\\
 \P^n\ar@{.>}[rr]^{\chi}&              &\P^n&}
 \]

The divisor $\Pi_{\Lambda_1}(\Psi(Z))$ is an hypersurface of
degree, say $j$. The projection $\Pi_x$ is birational only if 
$$\mult_{\Pi_{\Lambda_1}(x)}\Pi_{\Lambda_1}(\Psi(Z))=j-1$$

On the other hand $x$ is a general point of $\Psi(H_Z)$ therefore this
forces
\begin{equation}
\label{eq:mult}
\mult_{\Pi_{\Lambda_1}(\Psi(H_Z))}\Pi_{\Lambda_1}(\Psi(Z))=j-1
\end{equation}

\begin{claim} $\Pi_{\Lambda_1}(\Psi(H_Z))$ is a non degenerate
  generically smooth codimension 2
  subvariety of $\P^{n+1}$
\end{claim}
\begin{proof}[Proof of the claim]
The variety $\Pi_{\Lambda_1}(\Psi(H_Z))$ is degenerate only if there
is an element in $\G_{d-1,n,l-h}$. To check this, by the main result in \cite{AH1}, it is enough to verify the following inequality
$$\bin{n+d-1}n<\bin{n+d+1}{n}-\bin{n+d}{n-1}-(n+1)$$

Moreover $\Pi_{\Lambda_1}\Psi(H_Z)$ is a general projection from $h-1$
general points of $\Psi(H_Z)\subset\P^{n+h}$. Therefore by the trisecant Lemma, see
for instance
\cite[Proposition 2.6]{CC}, the map is birational.
In particular I have $\dim\Pi_{\Lambda_1}(\Psi(H_Z)=n-1$.
\end{proof}

The claim says that the secant variety of $\Pi_{\Lambda_1}(\Psi(H_Z)$
 fills up $\P^{n+1}$. Therefore by equation (\ref{eq:mult}) I have
$$2(j-1)\leq j$$
That is $j=2$ and $\Pi_{\Lambda_1}(\Psi(Z)$ is smooth along
$\Pi_{\Lambda_1}(\Psi(H_Z))$. The general choice of the points allows
to conclude that 
$$\deg\Psi(Z)=2+(h-1)=h+1$$

To get a lower bound on the degree of $\Psi(Z)$ I
consider a general hyperplane 
$$A\subset\langle \Psi(C_Z)\rangle\iso\P^{h+1}$$
 containing
$\Psi(E_i)\cap \Psi(C_Z)$. Then for such an $A$, by equation  (\ref{eq:span}) and (\ref{eq:birat}) I have
$$A\cdot\Psi(C_Z)\geq 2^{n-1}+h+1-\bin{n+1}2$$ 
A simple computation gives, for $n\geq 5$,
$$2^{n-1}>\bin{n+1}2$$
therefore for $n\geq 5$ I conclude
$$\deg\Psi(C_Z)>h+1$$
This proves that $\chi$ is finite and not birational. Then, by the
usual degeneration argument,
the map given by $\G_{d+1,n,l}$ is not birational.

Note that for $n=3,4$ the statement is a consequence of \cite[Theorem 1]{Me}. 
\end{proof}

I am in the condition to prove the main  Theorem
\begin{proof}[Proof of Theorem \ref{th:main}] 
The existence of a unique decomposition translates, by Proposition
\ref{pro:bridge}, in the fact that $\G_{d,n,l}$, with
$l=\frac{\bin{n+d+1}{n}}{n+1}-1$, gives a birational
map to $\P^n$. For $n\leq 4$ the result is proved in \cite{Me}. For
$n\geq 5$ using Theorem \ref{th:bsl} I conclude that such a linear system cannot exist.
\end{proof}

\begin{remark} If the map induced by sections of
  $\G_{d,n,l}$ is birational then a general hypersurface of degree $d$, with
  $l$ ordinary nodes, \cite{Me}, is rational. This is quite against general
  expectations, at least for $d$ big enough. Unfortunately I do not know any direct method to prove
  the non rationality of these special hypersurfaces. One should
  confront this with Koll\'ar result, \cite{Ko}, on non rationality of
  very general smooth hypersurfaces of degree roughly bounded by $2/3$ the dimension.
\end{remark}

\end{document}